\documentclass{ifacconf}
\usepackage{url}

\usepackage{hyphenat}
\usepackage{graphicx}      
\usepackage{natbib}        
\usepackage{graphicx} 
\usepackage{color}
\usepackage{amsfonts,amssymb}
\usepackage{bm}
\usepackage{amsmath}
\usepackage{multirow}
\usepackage{algorithmic}
\usepackage{graphicx}
\usepackage{wrapfig}
\usepackage{subcaption}
\usepackage{float}
\usepackage{textcomp}
\usepackage{graphicx} 
\usepackage{epstopdf}
\usepackage{array}
\newtheorem{mythm}{Theorem}
\newtheorem{mydef}{Definition}

\newtheorem{proof}{Proof}
\newtheorem{assumption}{Assumption}
\newtheorem{lemma}{Lemma}
\newtheorem{remark}{Remark}
\newtheorem{example}{Example}

\usepackage[title]{appendix}

\usepackage{array}
\usepackage{enumerate}
\usepackage{booktabs}
\usepackage{algorithm}
\usepackage{algorithmic}
\usepackage{setspace}
\usepackage{float}
\usepackage{cases}

\usepackage{mathrsfs}

\definecolor{kc}{rgb}{0.2,0.8,0.2}

\makeatletter\@mparswitchfalse\makeatother\normalmarginpar 
\usepackage[textwidth=0.8cm, textsize=scriptsize]{todonotes}

\begin{document}
\begin{frontmatter}

\title{ Game-Theoretic Learning-Based Mitigation of Insider Threats
\thanksref{footnoteinfo}} 

\thanks[footnoteinfo]{This work is supported in part by NSF under Grants CNS-2401007, CMMI-2348381, IIS-2415478, in part by MathWorks, and in part by EPSRC [grant number EP/X033546].}

\author[First]{Gehui Xu},  
\author[Second]{Kaiwen Chen}, 
\author[Second,Third,Fourth]{Thomas Parisini}, 
\author[First]{Andreas A. Malikopoulos}

\address[First]{School of Civil $\&$ Environmental Engineering, Cornell University, Ithaca, NY, USA}
\address[Second]{Department of Electrical and Electronic Engineering, Imperial College London, UK}
\address[Third]{Department of Electronic Systems, Aalborg University, Denmark}
\address[Fourth]{Department of Engineering and Architecture, University of Trieste, Italy}

\begin{abstract}                

An insider is defined as a team member who covertly deviates from the team’s optimal collaborative control strategy in pursuit of a private objective, while maintaining an outward appearance of cooperation. Such insider threats can severely undermine cooperative systems: subtle deviations may degrade collective performance, jeopardize mission success, and compromise operational safety.
This paper presents a comprehensive framework for identifying and mitigating insider threats in cooperative control settings. We introduce an insider-aware, game-theoretic formulation in which the insider’s hidden intention is parameterized, allowing the threat identification task to be reformulated as a parameter estimation problem. To address this challenge, we employ an online indirect dual adaptive control approach that simultaneously infers the insider’s control strategy and counteracts its negative influence. By injecting properly designed probing signals, the resulting mitigation policy asymptotically recovers the nominal optimal control law -- one that would be achieved under full knowledge of the insider’s objective.
Simulation results validate the effectiveness of the proposed identification–mitigation framework and illustrate its capability to preserve team performance even in the presence of covert adversarial behavior.

\end{abstract}

\begin{keyword}
 Insider threats, Game theory, Adaptive systems
\end{keyword}

\end{frontmatter}

\section{Introduction}
Ensuring safety and security in intelligent systems has received substantial 
attention, with extensive research on secure coordination \citep{farokhi2017private}, adversarial 
learning \citep{chen2024approaching}, and threat diagnosis mechanisms \citep{zhang2025threat} across 
domains such as intelligent transportation \citep{malikopoulos2021optimal}, 
power systems \citep{higgins2020stealthy}, and physical human--robot interaction 
\citep{sheng2025human}. While much of this work focuses on external adversaries, 
insider threats, where deceptive or opportunistic non-cooperative behaviors arise 
from within the team, have become increasingly prominent. According to a recent 
global report \citep{PonemonDTEX2025}, the average annual cost associated with 
insider-related incidents increased by nearly 50\% between 2019 and 2025. An insider is a trusted team member who possesses legitimate, often privileged, access to internal resources, where the team is understood as an organization whose members work collaboratively toward a shared objective~\citep{radner1962team,Malikopoulos2021}. Although insiders may outwardly appear to support collective tasks, they can covertly exploit their privileged position to manipulate coordination or deliberately disrupt operations for personal, financial, or other improper gain~\citep{nurse2014understanding,cappelli2012cert}. Such concealed behaviors leave other team members unaware of the insider’s true intentions, resulting in misinformed decisions, degraded collective performance, and, in many cases, compromised operational safety.

The resulting behaviors can severely degrade team performance and, in some cases, jeopardize team safety. For example, in a cooperative lane-change or lane-merge scenario~\citep{zhang2024stackelberg,rios2016survey}, a following vehicle may appear to decelerate to create a merging gap, yet covertly accelerate during the maneuver to sideswipe the leading vehicle, thereby shifting collision liability and facilitating insurance fraud. In human–robot collaborative tasks~\citep{sheng2025human,wang2022bounded}, a human operator may seem cooperative but intentionally underexert effort to conserve energy, shifting the workload to the robot and increasing the risk of task failure, such as dropping a shared object. In microgrids, a malicious insider may leak sensitive topology information to external adversaries, enabling false-data injection attacks that manipulate voltage or current measurements and potentially trigger instability or widespread outages \citep{Gonen2020FDI}.

Since an insider aims to avoid detection by behaving in ways that appear consistent with normal collaboration, such threats often go unnoticed until substantial damage has already occurred. Consequently, identifying insider threats and developing effective mitigation strategies have become critical challenges. Interactions between an insider and other team members are frequently modeled using game-theoretic frameworks~\citep{liu2020defense,liu2021flipit,hu2015dynamic,xu2024consistency,liu2008game}. However, most existing models assume that cooperative agents have full knowledge of the insider’s behavior—an assumption that rarely holds in practice. Furthermore, effective mitigation must occur concurrently with intention identification, rather than relying solely on delayed corrective actions taken after the insider has already achieved its objective.


The main objective of this paper is to address this gap by developing an integrated 
identification--mitigation framework for insider threats. We begin by formulating a 
two-player insider-aware game model built upon the nominal team-decision structure, 
and we parametrize the insider's hidden intention so that the threat-identification 
problem can be recast as one of parameter estimation. Building on this formulation, 
we design an online identification--mitigation scheme that concurrently infers the 
insider's intention from observed interactions while actively counteracting its 
adverse influence.
The resulting mitigation control strategy asymptotically coincides with the optimal one that assumes full knowledge of the insider’s intention.   Simulation results validate the effectiveness of the proposed approach. 


 \textit{Notation.} Let \(\mathbb{N}^+\) denote the set of positive integers and \(\mathbb{R}\) denote the set of real numbers.
We use \(\mathbb{R}^{n}\) (or \(\mathbb{R}^{m\times n},\; m,n\in\mathbb{N}^+\)) to denote the set of \(n\)-dimensional real column vectors (or real \(m\)-by-\(n\) matrices). 
Let \(I_n\) denote the \(n \times n\) identity matrix, $\otimes$ denote the Kronecker product, and \(\nabla f\) denote the gradient of a differentiable function \(f\).
For a signal \(x(t)\), we write \(x \in \mathcal{L}_\infty\) if it is bounded, \textit{i.e.}, 
\(\sup_{t\ge 0} \|x(t)\| < \infty\), and \(x \in \mathcal{L}_2\) if it is square-integrable, \textit{i.e.}, 
\(\int_{0}^{\infty} \|x(t)\|^2 \, dt < \infty\).

\section{A Two-player team game}

In this section, we formalize the collaborative setting introduced above by 
transitioning from the notion of a \emph{team} of cooperating members to an 
equivalent \emph{two-player game} representation. This game-theoretic formulation captures the strategic interaction between the decision maker  (representing the cooperative team behavior) and the potential insider.  We first model the nominal collaborative scenario in the absence of insider 
threats and then introduce the corresponding insider-threat formulation, in which 
the insider deviates from the team objective. Finally, we provide two illustrative 
examples that highlight how these formulations apply in practical settings.

\subsection{Nominal Scenario}

Consider a  two-player dynamic system described by
\begin{equation}\label{dynamics}
\dot{x} =f(x)+g_1(x)u_1+g_2(x)u_2 \, ,
\end{equation}
where $x(t) \in \mathbb{R}^n$, $n\in\mathbb{N}^+$, is the state of the system, 
$f(x): \mathbb{R}^n \to \mathbb{R}^n$, 
$g_1(x) \in \mathbb{R}^{n\times m}$ and 
$g_2(x) \in \mathbb{R}^{n\times m}$, $m \in\mathbb{N}^+ $, are smooth mappings, 
and \(u_1 \in U_1\subseteq\mathbb{R}^{m}\), \(u_2 \in U_2\subseteq\mathbb{R}^{m}\) are the control inputs of players 1 and 2, respectively, where $U_1$ and $U_2$ are given compact sets. Both players are assumed to have full observation of the system state and can apply control actions to influence its evolution. 

The nominal collaboration model is formulated as a team game in which two players jointly accomplish a task through cooperation and have the same cost function.
 Accordingly, we refer to player~1 as the decision maker (DM) and player~2 as the insider.
The common cost function for players to minimize is:
\begin{align}\label{cooperative}
\mathcal{C}\!=\! \int_0^\infty \!\!
&(x(t)\!-\!x^{r}_c)^\top\! Q_c(x(t)\!-\!x^{r}_c)\!+\! u_1(t)^\top \! R_{1} u_1(t)
\!\notag\\
&+\! u_2(t)^\top \!R_{2} u_2(t)\,dt,   
             \end{align}
             where $Q_c \succ 0$ and $R_{1},R_{2} \succ 0$ are weighting matrices, and $x_c^r \in \mathbb{R}^n$ denotes the desired reference state associated with the team objective. 

\begin{remark} Common cost functions such as~\eqref{cooperative} may represent a variety of practical scenarios. For instance, in a lane-change scenario, the leading vehicle and the following vehicle aim to maintain a desired safety distance and a preferred speed. Or, in a human–robot interaction scenario, a robot and a human jointly move a heavy object  to a desired location. Further details are provided in Section~\ref{sec:example}.
\end{remark}


The team outcome resulting from the joint decisions of all players is characterized by the team-optimal solution~\citep{radner1962team,zoppoli2020neural,Malikopoulos2021}. This solution corresponds to a strategy profile under which  no unilateral or joint deviation by players can yield improved collective performance~\citep{xu2025does,Malikopoulos2021}.
We consider team-optimal solutions under a dynamic closed-loop information structure, 
where the DM and the insider choose instantaneous control strategies 
$u_1(x)$ and $u_2(x)$ based on the observed state $x$.  To formalize the class of feasible feedback strategies, we define the admissible control spaces for the DM and insider as follows:
\[
\begin{aligned}
\mathcal{U}_1 &\!=\! \{u_1 \!\mid\! u_1\!:\!\mathbb{R}^n \!\rightarrow U_1,u_1(x)\ \text{is continuous in}\ x\},\\
\mathcal{U}_2 &\!=\! \{u_2 \!\mid\! u_2\!:\!\mathbb{R}^n \!\rightarrow U_2, u_2(x)\ \text{is continuous in $x$}\}.
\end{aligned}
\]
Next, we formalize the notion of optimal collaboration within the nominal team 
setting. 

\begin{mydef}
A pair $(u_1^*,u_2^*) \in \mathcal{U}_1 \times \mathcal{U}_2$ is called a 
team-optimal solution if
\[
\mathcal{C}\big(u_1^*,u_2^*;x_0\big)
\;\le\; \mathcal{C}\big(u_1,u_2;x_0\big),
\quad \forall\, u_1 \in \mathcal{U}_1,\; u_2 \in \mathcal{U}_2.
\]
\end{mydef}

\subsection{Insider Threat Scenario}\label{sec:insider}


In practice, an insider may appear to contribute to the team's task while secretly optimizing their own objective, with the DM unaware of this intention. Instead of following the original team cost, the insider optimizes an alternative objective that blends the team-oriented component with a private objective reflecting selfish or malicious intention:
\begin{align}\label{insider_obj}
\mathcal{C}^{\text{adv}}_2\!=\!&\int_0^\infty \!\!\!(x(t)-x_a^r)^\top Q_a (x(t)-x_a^r)
\!+ \!u_2(t)^\top\! \tilde R_{2} u_2(t)\!
\\
&+
\rho(u_2(t)\!-\!u_2^{*}(t))^\top\! (u_2(t)\!-\!u_2^{*}(t)) dt \notag
\end{align}
where $Q_a \succ 0$, $\tilde R_{2} \succ 0$ are insider-specific weighting matrices, and $x_a^r\in \mathbb{R}^n$ denotes the insider’s preferred reference state, which may differ from $x^r_c$. {These parameters, together with \(\rho\), are unknown to the DM.}
The last term, referred to as the {disciplinary risk}, 
penalizes the deviations of the insider’s control strategy from the nominal team strategy \(u_2^{*}(t)\), 
with coefficient \(\rho>0\).  {A larger value of $\rho$ corresponds to more cautious and concealed behavioral changes, whereas a smaller $\rho$ reflects more aggressive and overt adversarial actions.}

Accordingly, the insider knows that the DM is unaware of its true intention and selects its best response to the DM's nominal strategy $u_1^*$:
\[
u_2^\diamond=u_2^\diamond(u_1^*)
\;\in\;
\arg\min_{u_2 \in \mathcal{U}_2}
\mathcal{C}_2^{\mathrm{adv}}(u_1^*,u_2;x_0).
\]
Without awareness of the insider threat, the DM cannot prevent manipulation or disruption of the coordination process and therefore executes \(u_1^*\) under the false assumption of cooperation. This strategic asymmetry can lead to significant performance degradation or even more severe consequences.
To address such threats, the DM should infer the insider’s actual behavior and accordingly devise an effective mitigation control strategy. Therefore, the question of interest is:
how can the DM identify and mitigate such insider threats?

Now, we provide two examples that illustrate the specific forms and applicability of the nominal and insider-threat models.

\subsection{Examples}\label{sec:example}

\begin{example}[Vehicle Lane  Change]
Consider a lane\hyp changing scenario where the leading vehicle attempts to merge in front of the following vehicle \citep{zhang2024stackelberg,falsone2022lane}.
The following vehicle, acting as an insider, pretends to slow down to allow the leading vehicle to merge. However, it then accelerates while the leading vehicle is changing lanes, attempts to lightly collide with it, making the leading vehicle appear at fault and enabling an insurance claim. 

For each vehicle, consider the longitudinal kinematics
\[
\dot p_i = v_i,\qquad \dot v_i = u_i,\qquad i\in\{1,2\},
\]
where \(p_i\) and \(v_i\) denote position and velocity, and \(u_i\) is the
 acceleration strategy.
Define the system state as
$
x = [p_1-p_2,v_1,v_2]^\top \in \mathbb{R}^3,
$
which evolves according to
\[
\dot x \!=\! A x + B_1 u_1 + B_2 u_2,
\,
A \!=\!
\begin{bmatrix}
 0 & 1 & -1\\
 0 & 0 & 0\\
 0 & 0 & 0\\

\end{bmatrix},
B_1 \!=\!
\begin{bmatrix}
 0\\ 1\\ 0
\end{bmatrix},
B_2 \!=\!
\begin{bmatrix}
 0\\ 0\\ 1
\end{bmatrix}.
\]
In the absence of insider threats, the leading vehicle and the following vehicle aim to  maintain a desired safety distance  \(\pi>0\) and a desired
velocity \(v_c>0\), while minimizing their control efforts during the lane-change maneuver~\citep{zhang2024stackelberg}.
The common cost function is defined as
\begin{align}\label{team_transportation}
\mathcal{C}
\!=\!\!\! \int_{0}^{\infty}\!  (x(t)\!-\!x^r_{c})^\top Q_c(x(t)\!-\!x^r_{c}) \!+\!r_1u_1^2(t)+\!r_2u_2^2(t) dt,
\end{align}

where $Q_c =
\mathrm{diag}(q_1, q_2, q_2)$, $x^r_{c}= \begin{bmatrix} \pi,~ {v}_c,~ {v}_c\end{bmatrix}^{\!\top}$, and \(q_1> 0\), \(q_{2}> 0\), \(r_1 > 0\), \(r_2 > 0\) are weighting
parameters.

However,  the following vehicle may pursue a self-interested objective, and its true cost function is given by
\begin{align}\label{traffic:insider}
\mathcal{C}_2^{\text{adv}}\!=\!\! \int_{0}^{\infty}\!  (x(t)\!-\!x^r_a)^\top Q_a(x(t)\!-\!x^r_a) \!+\!\tilde r_2u_2^2\!+\!\rho(u_2-u_2^{*})^2 dt,
\end{align}
where $
Q_a =\mathrm{diag}(\tilde{q}_1, \tilde{q}_1, \tilde{q}_2 )
$, and $ x^{r}_a= \begin{bmatrix} \tilde{\pi},~ {v}_{a},~ {v}_{a}\end{bmatrix}^{\top} 
$.

Here, \(\tilde{\pi}\geq 0\) denotes the actual spacing that the following vehicle intends 
to maintain with the leading vehicle. 
When the following vehicle attempts to deliberately cause a collision, we have
\(\tilde{\pi}=0\). 
The term \(u_2^{*}\) represents the optimal control strategy 
that vehicle~2 is supposed to follow in the nominal case.
The parameters \(\tilde{q}_1\ge0\), \(\tilde{q}_{2}\ge0\), and \(\tilde{r}_2>0\) 
are the corresponding weighting coefficients, 
and \({v}_a>0\) denotes the target velocity that the following vehicle actually prefers to achieve.
\end{example}

\begin{example}[Human-Robot Interaction]
Consider a collaborative task in which a human and a robot jointly move a heavy object (e.g., a couch) to a desired position \citep{sheng2025human,wang2022bounded}. The human appears to cooperate but intentionally exerts less effort to save energy, thereby leading to an insider-threat scenario. 

In the nominal setting, the cost function jointly minimized by the robot (DM) and the human (insider) consists of a payoff term related to goal tracking and two penalty terms on their control efforts:
$$
\mathcal{C}\!=\!\int_{0}^{\infty}\! (x(t)-x^{r}_c)^\top Q_c(x(t)-x^{r}_c) +\omega_r\|u_1(t)\|^2+\omega_h\|u_2(t)\|^2.
$$
The system state is  $
x = [x_1, x_2] \in \mathbb{R}^6
$,  where \( x_1=[p_x, p_y, \sigma]^\top \in \mathbb{R}^3\) denotes the position of the object's center of mass and orientation, and 
\(x_2=[v_x, v_y, \dot{\sigma}]^\top \in \mathbb{R}^3\) represents the translational and angular velocities. Following \cite{lawitzky2010load,mortl2012role}, the system evolves as
$$
\dot{x}= f(x)+g_1(x)u_1+g_2(x)u_2,
$$
where
\[
f(x)\!=\!\!
\begin{bmatrix}
x_2\\
-M_o^{-1} f_o(x_1,x_2)
\end{bmatrix},
g_i(x)\! =\!\!
\begin{bmatrix}
0_{3\times 2} \\
M_o^{-1} G_i(x)
\end{bmatrix}, i\!\in\!\{1,2\},
\]
\(M_o\in\mathbb{R}^{3\times3}\) is the positive definite inertia matrix of the planar rigid body, 
and \(G_i(x)\in \mathbb{R}^{3\times2}\) denotes the partial grasp matrix associated with the \(i\)-th player, relating its applied force to the resultant object wrench.  
The term \(f_o(x_1,x_2)\in\mathbb{R}^3\) models passive effects such as ground friction and environmental interactions.  
Since the manipulated object is bulky and low sensitivity to grasp-induced torques, the interaction inputs are modeled as planar forces \(u_1,u_2\in\mathbb{R}^2\).

The parameter \(\omega_l>0\), \(l\in\{r,h\}\), characterizes each player's attitude toward effort.
A smaller \(\omega_l\) indicates a more active collaborator who is willing to exert stronger control inputs, 
whereas a larger \(\omega_l\) corresponds to a more cautious one. 
Since the human seeks to reduce exertion, it may intentionally adopt a larger effort-weight parameter \(\tilde{\omega}_h\) while attempting to disguise the reduced effort as normal cooperative behavior, resulting in the following objective function:
$$
\begin{aligned}
& \mathcal{C}^{\text{adv}}_2= \int_{0}^{\infty}  (x(t)-x^{r}_c)^\top Q_a(x(t)-x^{r}_c) \\
&\quad\quad\quad\quad+ \tilde{\omega}_h \|u_2(t)\|^2
  +\rho \big\|u_2(t)-u_2^{*}(t)\big\|^2 \, dt.
  \end{aligned}
$$
\end{example}

\section{Learning-Based Mitigation}

In this section, we focus on the identification and mitigation of insider threats.
We first provide an online learning scheme to estimate the unknown parameters associated with the insider’s intention.
Then, based on the estimated parameters, a corresponding mitigation control strategy for the DM is developed.

\subsection{Identification of insider threat}

When the insider behaves selfishly, its intention is implicitly encoded in the 
unknown parameters of its adversarial objective \(\mathcal{C}^{\mathrm{adv}}_2\), 
which influence the system only through the resulting control strategy 
\(u_2^{\diamond}\). To make this strategy identifiable, we assume that the DM 
knows the functional structure of \(u_2^{\diamond}(\cdot)\) but not the specific 
values of the underlying parameters. By observing the system trajectory \(x(t)\) 
and incorporating its own  actions \(u_1\), the DM  seeks to infer the 
parameterized structure of \(u_2^{\diamond}\) and consequently design an appropriate 
mitigation strategy.

The key idea enabling such identification 
is to construct a {linear parametric model}, also known as a linear regression model, such that 
the unknown parameter vector appears linearly in an equation where all other signals and parameters are known.
In what follows, we present the details in a linear–quadratic (LQ) setting.

For LQ team games, it is well known that the problem admits a unique optimal 
feedback solution provided that the pair $(A, [B_1, B_2])$ is stabilizable\footnote{Together 
with $R_{1}\succ0$, $R_{2}\succ0$, and $Q\succ0$, this condition is necessary and 
sufficient for the existence of a unique solution to the team problem. We 
additionally assume that the control constraints are inactive.}. Throughout this 
work, we assume that this stabilizability condition holds.
The feedback law is given by
\begin{equation}\label{eq:nominal_feedback}
    {u}_i^* = -{K}_i^* x - k_i^*, 
    \quad i \in \{1,2\},
\end{equation}
where \(K_i^* = R_i B_i^\top P^*\) denotes the optimal state-feedback gain of player \(i\), and 
\(
k_i^* = -K_i^* x_c^r.
\)
The matrix \(P^*\) is the positive definite solution to the algebraic Riccati equation
$
    A^\top P^* + P^*A + Q_c - P^* B R^{-1} B^\top P^* = 0,
$
where \(B = [B_1 \; B_2]\) and \(R = \operatorname{diag}(R_1, R_2)\).


By substituting \(u_2^*\) into the insider’s true objective function \(\mathcal{C}_2^{\text{adv}}\), 
and noting that \(Q_a \succ 0\), \(\tilde{R}_2 \succ 0\), and \(\rho>0\), 
the insider’s genuine optimal response to the DM’s action \(u_1^*\) 
retains a linear structure \citep{anderson2007optimal}, provided the pair \((A - B_1 K_1^*,\, B_2)\) is stabilizable, i.e,
\begin{equation}\label{eq:true_feedback}
    {u}_2^\diamond = -{K}_2^\diamond x - k_2^\diamond,
\end{equation}
where $K_2^\diamond
= (\tilde{R}_2 + \rho I)^{-1}\big(B_2^\top P^\diamond + \rho K_2^*\big)$, ${k}_2^\diamond=-K_2^\diamond x_a^r$ and  $P^\diamond$ is the solution of the following algebraic Riccati equation equation 
\begin{align}\label{eq:riccati_adversarial}
&  \quad (A-B_1K_1^*)^\top P^\diamond + P^\diamond  (A-B_1K_1^*) + Q_a+\rho K_2^{*\top}  K_2^*\notag\\
   & - \!(P^\diamond B_2 \!+ \!\rho K_2^{*\top} )(\tilde{R}_2 \!+\! \rho I)^{-1}(B_2^\top P^\diamond \!+\! \rho K_2^*)      \!=\! 0.
    \end{align}
Here, the reference state $x_a^r$ is designed such that the steady-state bias 
$(A-B_1K_1^*)x_a^r -B_1k_1^* $ in the closed-loop system $\dot{x}=(A-B_1K_1^*)x+B_2u_2-B_1k_1^*$ vanishes.

Substituting~\eqref{eq:true_feedback} into the system dynamics yields
\begin{equation}\label{eq:dynamics_true}
    \dot{x} = A x + B_1 u_1 - B_2 {K}_2^\diamond x - B_2 k_2^\diamond.
\end{equation}
From the DM’s perspective, the terms  $-B_2 {K}_2^\diamond$ and  $-B_2 k_2^\diamond$ are unknown and should be inferred. 


Rearranging terms in \eqref{eq:dynamics_true} gives
\[
\dot{x} - A x - B_1 u_1 = -B_2 {K}^\diamond_2 x- B_2 k_2^\diamond.
\]

Since the state derivative $\dot{x}$ is typically unavailable in practice, we apply a first-order low-pass filter $\frac{1}{s + \lambda}[\cdot]$ with $\lambda>0$ to both sides to obtain a realizable filtered representation \citep{ioannou1996robust}:
\[
\frac{1}{s+\lambda} [\dot{x} - A x - B_1 u_1]
= \frac{1}{s+\lambda} [-B_2 {K}^\diamond_2 x- B_2 k_2^\diamond].
\]

Then we have
\[
\xi_1:=\frac{1}{s+\lambda} [\dot{x}] = \frac{s}{s+\lambda} [{x}] ,
\quad
 \xi_2:=\frac{1}{s+\lambda} [-A x - B_1 u_1^*],
\]
$$ \eta_1 :=\frac{1}{s+\lambda}[x] ,
\quad \eta_2 := \frac{1}{s+\lambda}[1].
$$
Moreover, define the augmented regressor and parameter vector as
\[
\Theta^\star\! :=\!
\begin{bmatrix}
-B_2 K_2^\diamond \; & -B_2 k_2^\diamond
\end{bmatrix}
\in \mathbb{R}^{n\times (n+1)},
\,
\phi(t) \!:=\! 
\begin{bmatrix}
\eta_1 \\ \eta_2
\end{bmatrix}\!\in \mathbb{R}^{ (n+1)}.
\]
Then the unknown component can be written as
\[
-\,B_2 K_2^\diamond x - B_2 k_2^\diamond = \Theta^\star\, \phi(t).
\]
For parameter estimation, we vectorize $\Theta$ as
\[
\theta^\star \!:=\! \operatorname{vec}(\Theta^\star)\in \mathbb{R}^{ n(n+1)},
\;
\Phi(t) \!:=\! I_{n} \otimes \phi(t)^\top\!\in \mathbb{R}^{n\times n(n+1)},
\]
so that the filtered regression model becomes
\[
z(t):=\xi_1+\xi_2 = \Phi(t)\theta^\star,
\]
where $z(t)$ is the measured output signal, $\Phi(t)$ is the known regressor matrix, and $\theta^\star$ is the parameter vector unknown to the DM. 

Next, we illustrate how these unknown parameters can be estimated using an online 
parameter identification method. The central idea is to compare the measured system 
response \(z(t)\) with the output of a parameterized model \(\hat{z}(\theta,t)\) that 
shares the same structural form as the true plant dynamics. The parameter estimate 
\(\theta(t)\) is updated continuously so that \(\hat{z}(\theta,t)\) progressively 
matches the observed signal \(z(t)\) as time evolves. Under appropriate excitation 
conditions, the convergence of \(\hat{z}(\theta,t)\) toward \(z(t)\) guarantees that 
\(\theta(t)\) approaches the true parameter vector \(\theta^\star\) of the system.

Specifically, let $\theta(t)$ denote the estimate of the true parameter 
$\theta^\star$, with $\hat{z}(t)=\Phi(t)\theta(t)$ representing the predicted 
output. In contrast to classical gradient-based update laws that rely on a static 
adaptation gain, we adopt an identifier with a dynamic adaptation gain (DAG) 
structure~\citep{chen2024continuous,zhang2024dynamiccdc} to facilitate the tuning of the learning process, which is 
described by:
\begin{subequations}\label{eq:DAG_id}
\begin{align}
\epsilon&=\frac{z-\Phi\theta}{m^2},\\
    \dot{\xi} &= F\xi + G\epsilon,  \\
    \eta &= H\xi + \Gamma \epsilon,  \\
    \dot{{\theta}} &= \Phi^\top\eta, 
\end{align}
\end{subequations}
Here, $\epsilon(t)\in\mathbb{R}^{n}$  denotes the normalized estimation error based on $\theta(t)$, 
$m^2 = 1 + n_s^2$ with $n_s$ denoting the normalizing signal so that $\Phi/{m} \in \mathcal{L}_\infty$. Typical choices for $n_s$ include $n_s^2 = \phi^\top\phi$,  $n_s^2 = \operatorname{trace}(\Phi^\top \Phi)$ and 
$n_s^2 = \phi^\top P \phi$ with $P=P^\top \succ 0$. Moreover, $\xi(t)\in\mathbb{R}^{n_\xi}$ denotes the internal state of the DAG,
$\epsilon(t)\in\mathbb{R}^n$ is the input to the identifier, and
$\eta(t)\in\mathbb{R}^n$ is the filtered output. The matrices in~\eqref{eq:DAG_id} are selected such that 
$(F,G)$ is controllable and $(H,F)$ is observable, and the 
transfer function
$
    T_{(H\xi)\epsilon}(s) := H (sI-F)^{-1} G
$
is strictly positive real (SPR). In addition, $\Gamma = \Gamma^\top \succ 0$.

\begin{remark}

The identifier \eqref{eq:DAG_id} makes the parameter estimation error dynamics a negative-feedback loop of passive systems. The prediction error $\epsilon$ is first filtered by the DAG and then multiplied by the regressor, and the resulting signal is used to update the parameter estimate. Notice that~\eqref{eq:DAG_id} serves as a dynamic generalization of the 
static-gain update law $\dot{\theta} =\Phi^{\top} \Gamma \epsilon$, since\footnote{
Since both the regression model and the DAG rely on filtered signals, nonzero
initial filter conditions introduce exponentially decaying transient terms. As the
underlying filter dynamics are exponentially stable (a consequence of the SPR
condition), these transients vanish asymptotically and do not affect the stability
or convergence of the update law.}$\dot{\theta} = \Phi^{\top}T_{\eta\epsilon}(s)[\epsilon]$. 
Furthermore, the proposed identifier reduces to the classical static form if 
the internal state variables in~(\ref{eq:DAG_id}b) are removed and $H = 0$.
\end{remark}
\begin{remark}\label{convergence}
The DAG in ~\eqref{eq:DAG_id}  can be interpreted as a passive compensator that provides additional
design flexibility compared to classical update laws with static gains.  To illustrate this, consider the scalar time-invariant case where the DAG reduces to a first-order compensator of the form
$T_{\eta\epsilon}(s)=\gamma + \frac{\beta}{\alpha s + 1}$ with
$\alpha>0$, $\beta>0$, $\gamma>0$ \citep{chen2024continuous}. 
In this representation, the designer can shape the learning behavior by independently tuning the high-frequency gain $\gamma$, the DC gain $\beta + \gamma$, and the cut-off frequency (via the choice of $\alpha$), rather than relying on a single fixed adaptation gain.
\end{remark}

Next, we focus on the property of DAG and the convergence  of $\theta(t)$, in particular whether \(\theta(t)\) converges as \(t \to \infty\), 
 and if so, whether the
limit equals the true parameter $\theta^\star$. To establish this result, we impose
the following assumption.
\begin{assumption}[Persistence of Excitation]\label{assumption: PE}
The filtered regressor $\phi(t)$ is persistently exciting (PE), that
is, there exist $T_0>0$,   $\alpha_0>0$ and $\alpha_1>0$ such that, for all  $t\ge 0$,
\begin{equation*}\label{eq:PE}
\alpha_0I_{n+1}\preceq\int_{t}^{t+T_0}\!\phi(\tau)\phi(\tau)^\top d\tau \;\preceq\; \alpha_1I_{n+1}.
\end{equation*}
\end{assumption}
The property of $\phi(t)$ in Assumption
\ref{assumption: PE} is referred to as PE, which is crucial in many adaptive schemes. This condition implies that the regressor signal $\phi(t)$ is sufficiently rich 
to excite all modes of the system, thereby guaranteeing parameter convergence of the estimated parameters to their true values.  
An obvious concern is that if $x(t)$ converges to $0$, $\phi(t)$ converges to the constant vector, which may render Assumption~1 unachievable. One way to address this issue is to add a probing signal to the control input, which leads to the so-called dual control~\citep{wittenmark1995adaptive, bhasin2013novel}.

Furthermore,  according to Lemma~2 in~\citep{chen2024continuous}, the DAG described by~(\ref{eq:DAG_id}b)–(\ref{eq:DAG_id}c)  is both strictly passive and input strictly passive with respect to the input 
$\epsilon$ and the output $\eta$. 
This useful property, together with the  PE condition, implies the bounded internal signals, vanishing prediction error, and asymptotic convergence of ${\theta}(t)$ to $\theta^\star$  \citep[Proposition 2, Theorem 2]{chen2024continuous}, as stated below.


\begin{lemma}\label{l1}
Consider the update law \eqref{eq:DAG_id}. Then the following properties hold:
\begin{enumerate}[(i)]
    \item All signals within the identifier are bounded, \textit{i.e.}, they belong to $\mathcal{L}_\infty$.
    \item $\theta\in \mathcal{L}_\infty$, $\dot{\theta}\in \mathcal{L}_2\cap \mathcal{L}_\infty$, $\epsilon\in \mathcal{L}_2\cap \mathcal{L}_\infty$, $\epsilon n_s\in \mathcal{L}_2\cap \mathcal{L}_\infty$,
    \item Under Assumption~\ref{assumption: PE}, and provided that $n_s,\phi \in \mathcal{L}_\infty$, 
    the estimation error $\tilde{\theta}(t):= \theta(t)-\theta^\star $ converges to 0 exponentially. 
\end{enumerate}
\end{lemma}

Lemma~\ref{l1} establishes that the parameter estimates converge to their true 
values. Consequently, the mitigation strategy constructed from the estimated 
parameters asymptotically coincides with the control law that would be implemented 
if the true parameters were known, as further discussed in the sequel.

\subsection{Mitigation control strategy design}
Based on the estimated parameters, the DM continuously updates its control policy. This results in an online mitigation mechanism, where the identification and control are executed simultaneously. 

From the DM’s perspective, the system evolves as
\begin{align}\label{mitigate_dys}
\dot{x} = A x + B_1 u_1 +{\Theta}_1(t) x + {\Theta}_2(t),
\end{align}
where ${\Theta}_1(t)$ and ${\Theta}_2(t)$ denote the DM's  estimates of 
the insider’s feedback term $B_2 K_2^\diamond$ and bias term $B_2 k_2^\diamond$ at time $t$,  and together form the previously 
defined parameter vector \(\theta(t)\).
To counteract the estimated influence of the insider, the DM modifies the  original team objective accordingly, yielding:
$$
\mathcal{C}_1^{\mathrm{mit}} = \!\! \int_{0}^{\infty}\!  (x(t)-x_m^r)^\top Q_m(x(t)-x_m^r) +u_1(t)^\top\tilde R_1u_1(t),
$$
where $Q_m \succ 0$, $\tilde R_{1} \succ 0$ are mitigation-specific weighting matrices, and $x_m^r\in \mathbb{R}^n$ denotes the  desired state under mitigation. 

For instance, in the lane-change scenario (Example~1), when the following vehicle accelerates in an attempt to close the gap, the leading vehicle responds by increasing its speed to enlarge the inter-vehicle distance and avoid collision, thereby restoring the desired safety distance with a slight adjustment in speed. A similar pattern appears in the human–robot collaboration case (Example~2). When the human deliberately reduces effort to conserve energy, the robot compensates by exerting additional force to ensure successful task execution while adjusting its posture to prevent the object from being dropped.

We next characterize the optimal mitigation control strategy corresponding to the modified objective function. The reference state $x_m^r$ is designed so that the steady-state bias 
$(A+\Theta_1^\star)x_m^r + \Theta_2^\star$ vanishes. Thus, 
provided that the pair \((A + \Theta_1^*, B_1)\) is stabilizable, 
the DM’s genuine mitigation control strategy to the insider's action \(u_2^\diamond\) admits the linear feedback form 
\(u_1^{\theta^\star} = -K_1^{\theta^\star}x - k_1^{\theta^\star}\), where 
\(K_1^{\theta^\star} = \tilde{R}_1^{-1} B_1^\top P^{\theta^\star}\) and 
\(k_1^{\theta^\star} = -K_1^{\theta^\star} x_m^r\), and \(P^{\theta^\star}\) is the unique solution to the associated Riccati equation. 

Accordingly, when the insider parameters are unknown and estimated online, provided that the pair \((A+\Theta_1(t), B_1)\) is stabilizable \citep{ioannou1996robust}, the adaptive mitigation control strategy takes the form
\[
u_1^{\theta(t)}(x) = -K_1^{\theta(t)}x - k_1^{\theta(t)},
\]
where \(K_1^{\theta(t)} = \tilde{R}_1^{-1}B_1^\top P^{\theta(t)}\), 
\(k_1^{\theta(t)} = -K_1^{\theta(t)} x_m^r\),
and {\(P^{\theta(t)}\) is the unique solution to the Riccati equation $(A+\Theta_1(t))^\top P^{\theta(t)}+P^{\theta(t)}(A+\Theta_1(t))-P^{\theta(t)}B_1\tilde{R}_1^{-1}B_1^\top P^{\theta(t)}+Q_m=0$. 
}





The remaining question is whether \(u_1^{\theta(t)}\) can achieve the same performance as \(u_1^{\theta^\star}\) as \(\theta(t) \to \theta^\star\).
 We now present the main result. The proof is given in Appendix. \ref{proof_thm1}.


\begin{mythm}\label{thm1}


Consider the closed-loop system~\eqref{mitigate_dys} operating under the mitigation 
control strategy $u_1^{\theta}(x(t))$. Then all trajectories of ~\eqref{mitigate_dys}
are bounded, and
the state $x(t)$ converges to the mitigation reference $x_m^r$ for any initial condition $x(0)\in\mathbb{R}^n$. Moreover, under 
Assumption~\ref{assumption: PE}, the mitigation strategy $u_1^{\theta}(x(t))$ 
converges to the optimal  strategy $u_1^{\theta^\star}(x(t))$ for all $t$.

\end{mythm}

The above identification–based mitigation approach follows from a certainty equivalence principle. The idea behind it is that as the parameter $\theta$ converges to the true one $\theta^\star$,  the performance of the adaptive controller $u_1^{\theta}$ tends to that achieved by $u_1^{\theta^\star}$ in the case of known parameters. This enables the DM to reliably infer and counteract the insider's behavior. In practice, the convergence speed of the parameter estimates can be further improved by adjusting the designed parameters in \eqref{eq:DAG_id} as mentioned in Remark \ref{convergence}.  This allows the DM to react more rapidly, thereby reducing the risk associated with the delayed mitigation of insider manipulation.




\section{SIMULATIONS STUDIES}

In this section, we validate the performance of the identification-mitigation scheme using the lane-change scenario introduced in Example 1. We first demonstrate the estimation accuracy of the inverse-learning procedure and then evaluate the performance of the resulting mitigation control strategy.

Consider a highway scenario where the leading vehicle aims to maintain a safe distance of 
\(\pi = 73\,\mathrm{m}\) \citep{highwaycode126} from the following vehicle during the lane-change maneuver, 
while the insider’s desired inter-vehicular distance is \(\tilde{\pi} = 0\,\mathrm{m}\).
This parameter $\tilde{\pi}$, together with the remaining parameters in the cost function~\eqref{traffic:insider}, 
is unknown to the DM and implicitly embedded in the insider’s true feedback policy 
\((K_2^{\diamond}, k_2^{\diamond})\). 

Since \(B_2=\begin{bmatrix}0,0,1\end{bmatrix}^\top\) only affects the last state component, all other rows of \(\Theta^\star\! :=\!
\begin{bmatrix}
-B_2 K_2^\diamond \!\quad\! & -B_2 k_2^\diamond
\end{bmatrix}\) are zero, and only the last row requires estimation. Therefore, we employ the reduced regressor $\tilde{\Phi}(t)
= 
\begin{bmatrix}
0 & 0 & 1
\end{bmatrix}^{\!\top}
\phi^\top(t)$
which does not affect the theoretical estimation guarantees.
Moreover, to satisfy the PE condition during identification,  a combination of small sinusoidal signals is injected into the leading vehicle’s control input to ensure sufficient information richness for learning.
The initial estimate $
\theta(0)=
\begin{bmatrix}
-\,B_2K_2^{*} \quad -\,B_2k_2^{*}
\end{bmatrix}
$ is chosen according to the DM’s nominal perception, \textit{i.e.}, assuming the insider follows the team control strategy.


 \begin{figure}[H]

    \begin{subfigure}[t]{0.98\linewidth}
		\centering
		\hspace{0.48cm}\includegraphics[width=8.1cm]{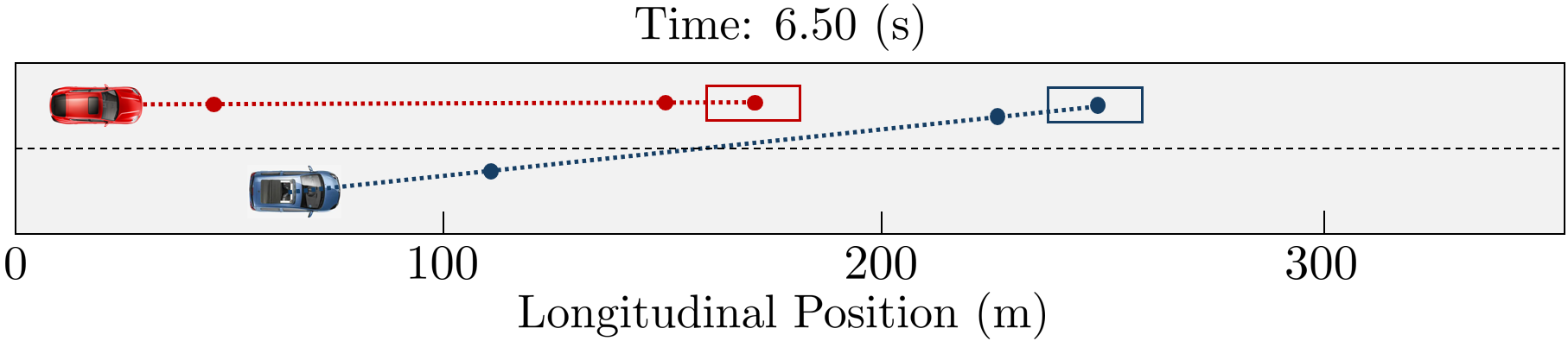}
	\end{subfigure}%
    \vspace{-0.02cm}
    \\
    	\begin{subfigure}[t]{0.98\linewidth}
		\centering
		\includegraphics[width=8.8cm]{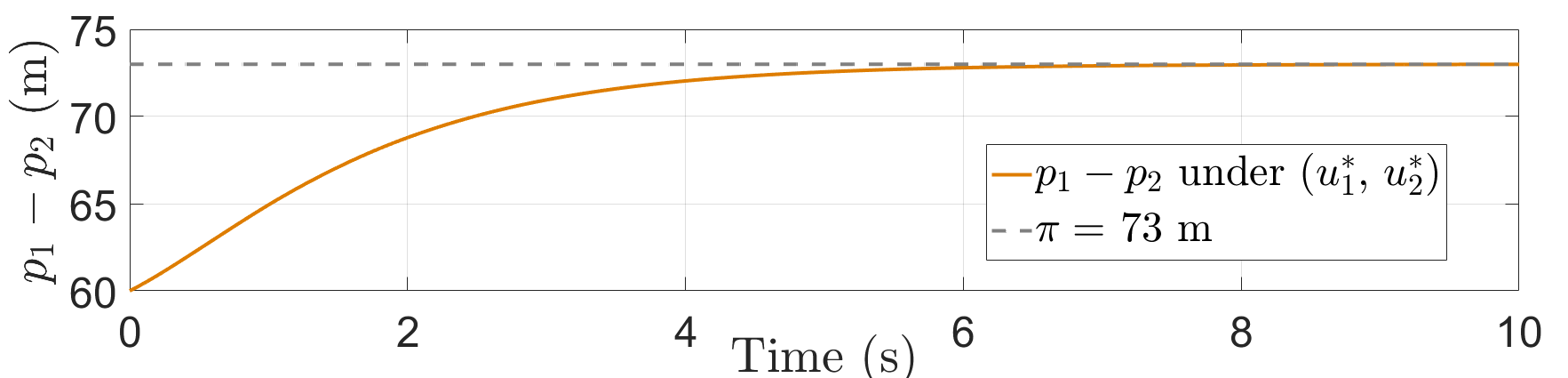}
	\end{subfigure}
	\centering
	\vspace{-0.2cm}
		\caption{The vehicles’ trajectory under team strategies.} 
	\label{fig22}
	\vspace{-0.2cm}
\end{figure}

\begin{figure}[H]

    \begin{subfigure}[t]{0.98\linewidth}
		\centering
		\hspace{0.48cm}\includegraphics[width=8.1cm]{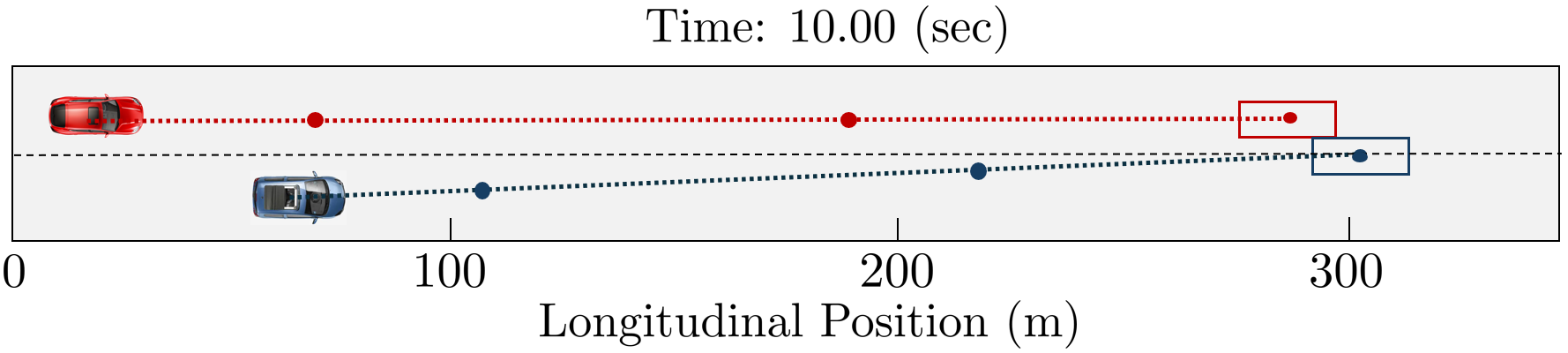}
	\end{subfigure}%
    \vspace{-0.02cm}
    \\
    	\begin{subfigure}[t]{0.98\linewidth}
		\centering
		\includegraphics[width=8.8cm]{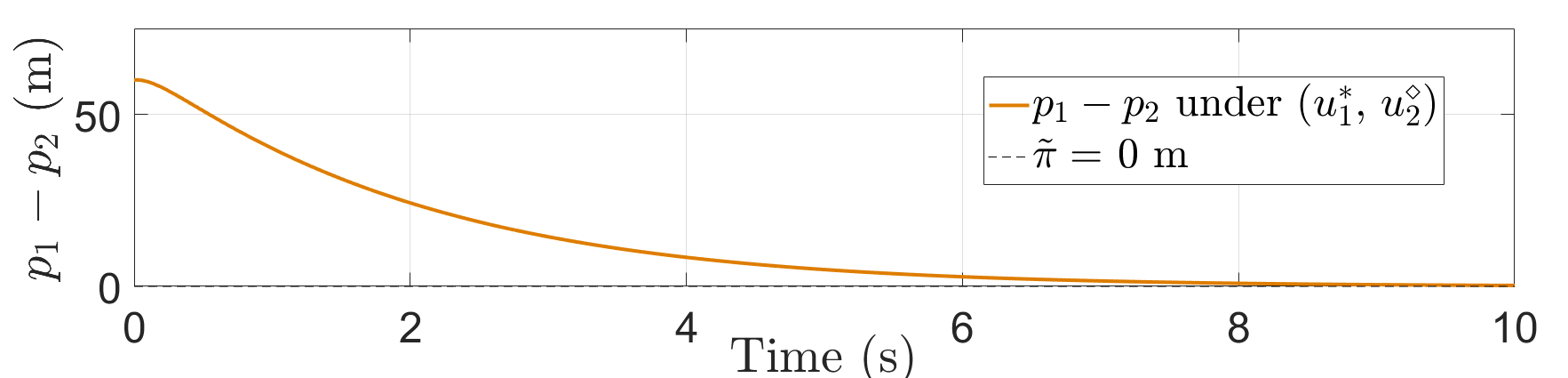}
	\end{subfigure}
	\centering
	\vspace{-0.2cm}
	\caption{The vehicles’ trajectory under insider threat.} 
	\label{fig14}
	\vspace{-0.2cm}
\end{figure}

Fig.~\ref{fig22} and Fig.~\ref{fig14} respectively illustrate the relative position trajectories in the nominal and insider-threat scenarios, with intermediate markers added at $t=2$~s and $t=6$~s. The former corresponds to a successful lane change, whereas the latter results in a sideswipe collision caused by the insider.

\begin{figure}[H]
	\centering	
			\includegraphics[width=0.98\columnwidth]{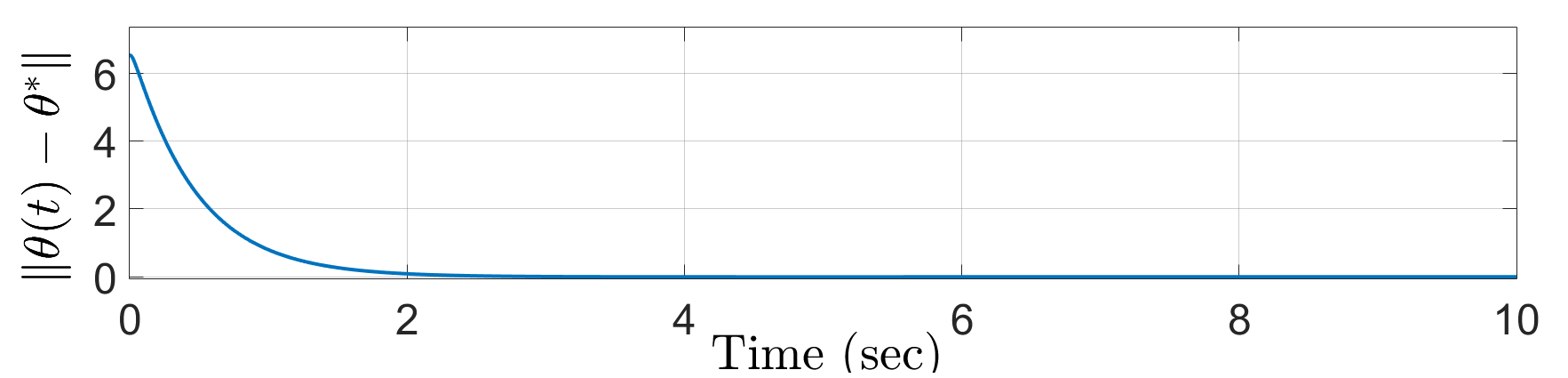}
\\
\vspace{-0.2cm}
	\caption{Insider intention learning performance} 
	\label{fig1}
\end{figure}

Fig.~\ref{fig1} shows the estimation error trajectory 
$\|\theta(t)-\theta^\star\|$ in the identification–mitigation scenario. The error converges to zero, demonstrating that the proposed learning scheme successfully identifies the insider’s true intention.



\begin{figure}[H]

    \begin{subfigure}[t]{0.98\linewidth}
		\centering
		\hspace{0.48cm}\includegraphics[width=8.1cm]{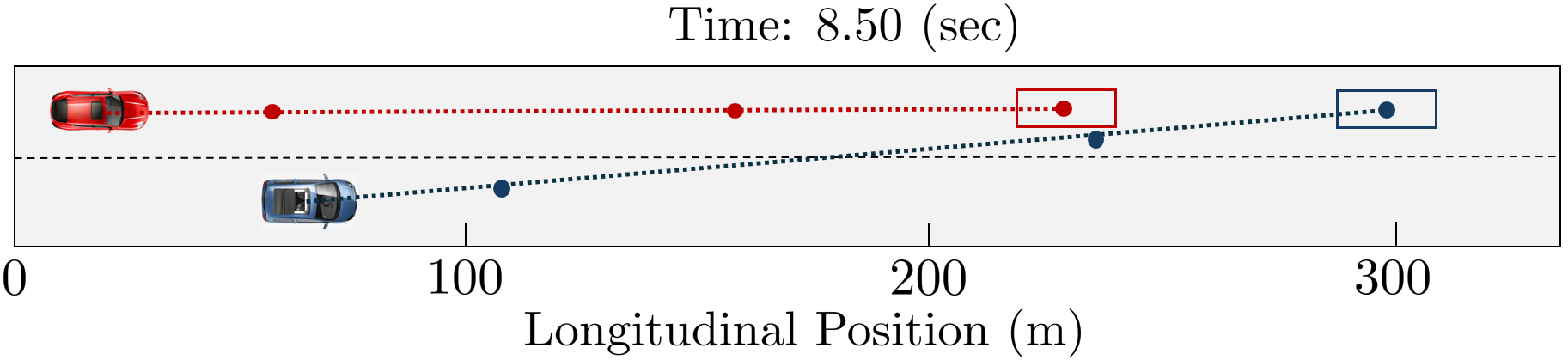}
	\end{subfigure}%
    \vspace{-0.02cm}
    \\
    	\begin{subfigure}[t]{0.98\linewidth}
		\centering
		\includegraphics[width=8.8cm]{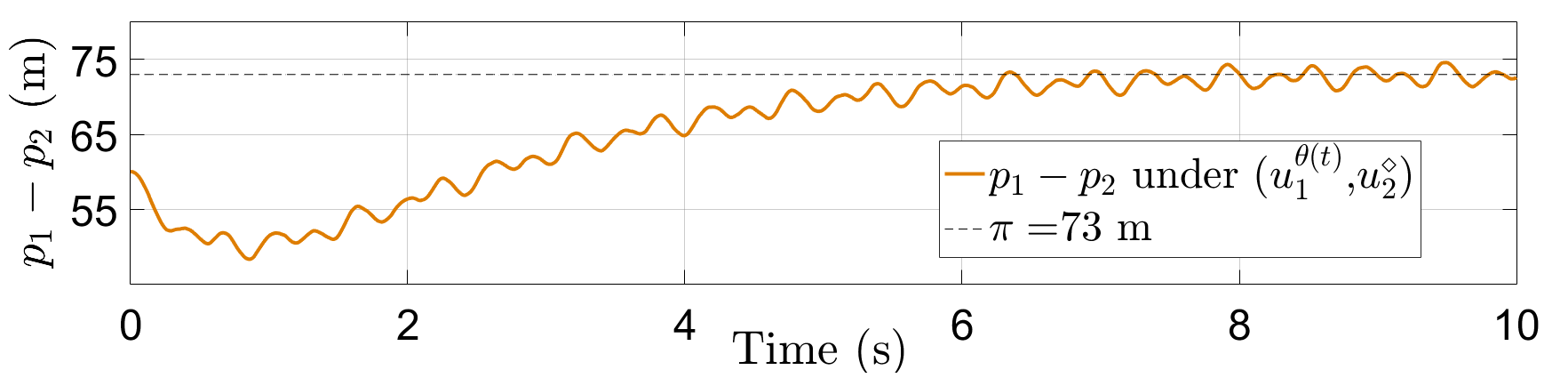}
	\end{subfigure}
	\centering
	\vspace{-0.2cm}
	\caption{The vehicles’ trajectory during the mitigation process.} 
	\label{fig3}
	\vspace{-0.2cm}
\end{figure}

Based on the learned parameters, the DM then computes the mitigation control strategy to ensure a safe distance when completing the lane change.
As the estimated  parameters converge to their true values, the corresponding feedback gains 
$K_1^{\theta(t)}$ and $k_1^{\theta(t)}$ also converge to their true counterparts 
$K_1^{\theta^\star}$ and $k_1^{\theta^\star}$, respectively.
Fig. \ref{fig3} presents the resulting trajectories of relative distance under the estimated mitigation control $u_1^{\theta}$. 
The leading vehicle gradually refines its control action and ultimately restores 
the safe \(73\,\mathrm{m}\) spacing. The oscillatory behaviors arise from the 
external sinusoidal excitation used for learning. This does not compromise the overall safety requirement as the amplitude of the oscillations is far smaller than the safety distance. On the other hand, the insider’s adversarial 
influence is successfully mitigated.

\begin{figure}[H]
	\centering	
			\includegraphics[width=1\columnwidth]{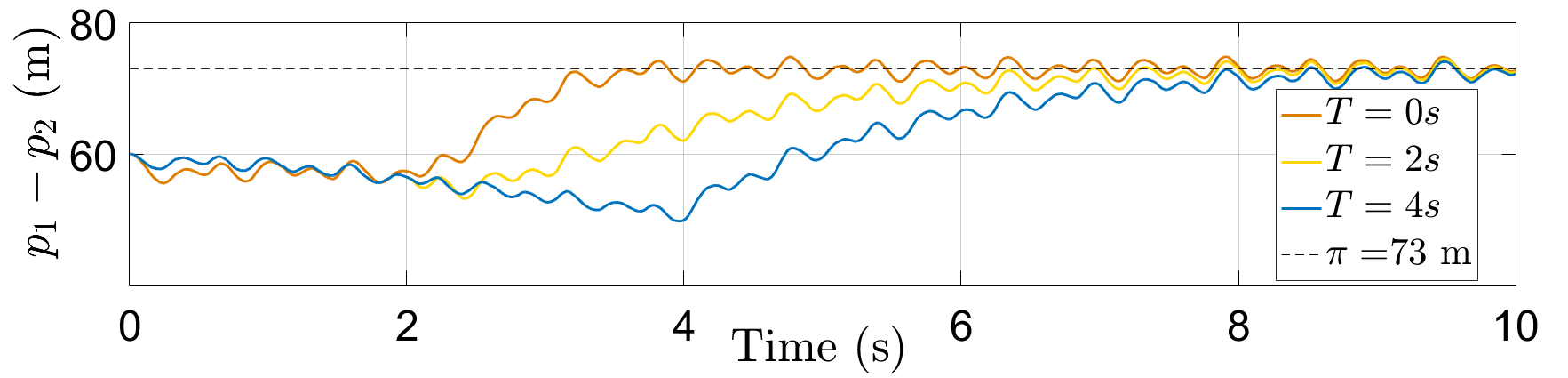}
\\
\vspace{-0.2cm}
	\caption{Impact of response timing on threat mitigation} 
	\label{fig16}
\end{figure}
Finally, we provide a discussion on the response timing in threat mitigation.
In practice, mitigation against an insider threat may not occur immediately, as threat identification often lags behind the onset of abnormal behavior. When the following vehicle slowly closes the gap, it may not trigger an immediate defensive reaction, causing the leading vehicle to continue the maneuver under the assumption of cooperative intention. As shown in the Fig.~\ref{fig16}, the leading vehicle triggers mitigation at 
$0$~s, $2$~s, and  $4$~s based on the detected trends in relative distance and relative velocity. It is evident that delaying the intervention increases the time required to re-establish the desired safety gap.  Furthermore, when the following vehicle exhibits increasingly aggressive intention, a feasible defensive response may involve aborting the lane change to ensure safety rather than attempting to recover spacing once risk becomes imminent.

\section{Concluding Remarks}
In this paper, we developed an identification-mitigation framework for addressing 
insider threats in cooperative systems. We formulated an insider-aware model within 
a game-theoretic setting and used a linear parameterization approach to recast the 
inference of hidden insider behavior as a parameter identification problem. Building 
on this formulation, we designed an online identification--mitigation scheme that 
simultaneously learns the insider's intention and counteracts its adverse influence. 
We further showed that the resulting mitigation control strategy asymptotically recovers 
the optimal control law that would be obtained under full knowledge of the insider's 
objective. In future work, we will extend the proposed framework beyond linear system 
settings to accommodate more general nonlinear and hybrid dynamics.


\bibliography{ifacconf}             

\appendix
\section{Proof of Theorem~1}\label{proof_thm1}    
\begin{proof}
 Consider the dynamics described by $\dot{x} = A x - B_1K_1^{\theta(t)}x  +{\Theta}_1(t) x + {\Theta}_2(t)$, and define  $e=x-x_m^r$, $\hat{A}(t)=A+{\Theta}_1(t)$, $A_{cl}(t)=\hat{A}(t)-B_1K_1^{\theta(t)}$, and $\omega(t)=(A+{\Theta}_1(t))x_m^r+{\Theta}_2(t)$. With some straightforward computation, one can see that the $e$-dynamics are described as follows. 
$$
\begin{aligned}
\dot{e}&=(A+{\Theta}_1(t))e+B_1u_1+(A+{\Theta}_1(t))x_m^r+{\Theta}_2(t)\\
&=A_{cl}(t)e+\omega(t).
\end{aligned}
$$
Since 
$(\hat{A},{B})$ is stabilizable at each time $t$, this implies that the unique solution 
$P^{\theta(t)} \succ 0$ of the Riccati equation exists and thus $P \in \mathcal{L}_\infty$, $K_1^{\theta(t)}\in \mathcal{L}_\infty$. 
Then one can establish that $A_{cl}(t)$ is a Hurwitz matrix at each 
frozen time $t$ as follows.  Computing the time derivative of $A_{cl}(t)$ yields
\begin{equation}
\| \dot{A}_{cl} \| 
\le 
\| \dot{\hat{A}} \| 
+ \|{B}_1\tilde{R}_1^{-1}{B}_1^\top\| \|\dot{P}^{\theta(t)}\|.
\end{equation}
Note that $\| \dot{\hat{A}} \| \in \mathcal{L}_2 $ since $\| \dot{\theta} \| \in \mathcal{L}_2 $ guaranteed by the update law.
Moreover, taking the time derivatives of both sides of the Riccati equation $\hat{A}^\top P^{\theta(t)}+P^{\theta(t)}\hat{A}-P^{\theta(t)}B_1\tilde{R}_1^{-1}B_1^\top P^{\theta(t)}+Q_m=0$, we obtain
\begin{equation}\label{Lyapunov}
\dot{P}^{\theta(t)}A_{cl} + A_{cl} ^\top \dot{P}^{\theta(t)} = -Q,
\end{equation}
where
\begin{equation}
Q = \dot{\hat{A}}^\top P^{\theta(t)} + P^{\theta(t)}\dot{\hat{A}}.
\end{equation}
Equation \eqref{Lyapunov}  is a Lyapunov equation, and its solution $\dot{P}^{\theta(t)}$ exists and  satisfies $\|\dot{P}^{\theta(t)}\| \le \nu \|Q(t)\|$ for any given $Q(t)$, where $\nu > 0$ is a constant.
Since $\|\dot{\hat{A}}(t)\| \in \mathcal{L}_2$ and $\hat{A}, P \in \mathcal{L}_\infty$, $K_1^{\theta(t)}\in \mathcal{L}_\infty$, we have $\|\dot{P}^{\theta(t)}\|\in \mathcal{L}_2$ and thus 
$\|\dot{A}_{cl}(t)\| \in \mathcal{L}_2$. 
Since \(A_{cl}(t)\) is Hurwitz for each frozen time \(t\) and $\|\dot{A}_{cl}(t)\| \in \mathcal{L}_2$, 
by Theorem~3.4.11 in~\cite{ioannou1996robust}, $x=0$ is a uniformly asymptotically stable (u.a.s.) equilibrium of the system $\dot x(t) = A(t)x(t)$. Moreover, 
by adding and subtracting $\Theta_1^\star$ and $\Theta_2^\star$ in $\omega(t)$, we obtain
$
\omega(t)=(A+\Theta_1^\star)x_m^r + \Theta_2^\star+(\Theta_1(t)-\Theta_1^\star)x_m^r+(\Theta_2(t)-\Theta_2^\star)
$, since $\theta\in\mathcal{L}_\infty $ and $\theta(t)\rightarrow \theta^{\star}$, we have
$\omega(t)\in\mathcal{L}_{\infty}$ with  $\omega(t)\rightarrow 0$.
Based on these and following 
the proof of \citep[Theorem~7.4.2]{ioannou1996robust}, it follows that 
\(e(t)\in\mathcal{L}_{\infty}\) and that \(e(t)\rightarrow 0\) as 
\(t\rightarrow\infty\).

On this basis, since the first-order filter \(1/(s+\lambda)\) is Bounded-Input Bounded-Output stable, we have  $n_s,\phi\in \mathcal{L}_\infty$.
Together with Assumption 1, the parameter estimation error
$\tilde{\theta}(t)$ decays exponentially according to Lemma~1~(iii).

Finally, define $\mathcal{S}(\theta):\theta\mapsto (K_1^{\theta},k_1^{\theta})$. When $\theta$ tends to $\theta^\star$, the continuity of $\mathcal{S}(\theta)$ implies $K_1^{{\theta}(t)}\rightarrow K_1^{\theta^\star}$ and $k_1^{{\theta}(t)}\rightarrow k_1^{\theta^\star}$.
Then for any state $x(t)$,
$u_1^{{\theta}(t)}(x(t))\to u_1^{\theta^\star}(x(t))$.  This completes the proof. \hfill$\square$
\end{proof}
\end{document}